\documentclass[letterpaper, 10 pt, conference]{ieeeconf}

\IEEEoverridecommandlockouts 
\usepackage{amsmath, amssymb, mathtools, mathrsfs, xcolor}
\usepackage[american]{circuitikz}
\usepackage{subcaption}
\usepackage{pgfplots}
\usepackage{tikz}

\ctikzset{bipoles/length=1.4cm}

\DeclareMathOperator{\im}{im}

\newcommand{\bbR}{\mathbb{R}}
\newcommand{\bbC}{\mathbb{C}}

\newcommand{\bbT}{\mathbb{T}}

\newcommand{\B}{\mathscr{B}}
\newcommand{\D}{\mathscr{D}}

\newcommand{\calA}{\mathcal{A}}
\newcommand{\calB}{\mathcal{B}}
\newcommand{\calC}{\mathcal{C}}
\newcommand{\calD}{\mathcal{D}}

\newcommand{\calF}{\mathcal{F}}
\newcommand{\calG}{\mathcal{G}}
\newcommand{\calH}{\mathcal{H}}

\newcommand{\calJ}{\mathcal{J}}

\newcommand{\calQ}{\mathcal{Q}}

\newcommand{\calR}{\mathcal{R}}
\newcommand{\calX}{\mathcal{X}}

\newcommand{\calT}{\mathcal{T}}
\newcommand{\calY}{\mathcal{Y}}
\newcommand{\calS}{\mathcal{S}}
\newcommand{\calU}{\mathcal{U}}

\newcommand{\calZ}{\mathcal{Z}}

\renewcommand{\epsilon}{ \varepsilon}

\newcommand{\bbm}{\begin{bmatrix*}}
\newcommand{\ebm}{\end{bmatrix*}}

\newcommand{\inv}{^{-1}}
\renewcommand{\t}{^\top}
\renewcommand{\a}{^\ast}

\newcommand{\norm}[1]{\lVert#1\rVert}
\newcommand{\inner}[2]{\langle#1, #2\rangle}

\newcommand{\half}{\frac{1}{2}}

\renewcommand{\d}{\textrm{d}}

\newcommand{\sys}{\Sigma}
\newcommand{\hatbar}[1]{\hat{\bar{#1}}}

\newcommand{\bbRp}{\bbR_{\geq 0}}

\newcommand{\intrs}{\int_{\bbR^s}}

\newcommand{\dual}{^{\d}}

\newtheorem{example}{Example}
\newtheorem{definition}{Definition}

\newtheorem{theorem}{Theorem}
\newtheorem{proposition}{Proposition}

\newtheorem{remark}{Remark}
\newtheorem{corollary}{Corollary}

\title{\LARGE \bfseries
	Passive and reciprocal linear time-and-space-invariant systems
}

\author{Brayan M. Shali$^1$, \quad Rodolphe Sepulchre$^{1,2}$
\thanks{$^\ast$The research leading to these results has received funding from the European Research Council under the Advanced ERC Grant Agreement SpikyControl n.101054323. Email: {\itshape brayan.shali@kuleuven.be}, {\itshape rodolphe.sepulchre@kuleuven.be}.}%
\thanks{$^1$Department of Electrical Engineering (ESAT), KU Leuven, KasteelPark Arenberg 10, B-3001 Leuven, Belgium.
}
\thanks{$^2$Department of Engineering, University of Cambridge, TrumpingtonStreet, Cambridge CB2 1PZ, United Kingdom.}%
}

\begin{document}
	
	\maketitle
	\thispagestyle{empty}
	\pagestyle{empty}
	
	\begin{abstract}
		Reciprocity is a fundamental symmetry property observed across many physical domains, including acoustics, elasticity, electromagnetics, and thermodynamics. In systems and control theory, it provides key insights into the internal structure of linear time-invariant (LTI) systems and is closely linked to properties such as passivity, relaxation, and time-reversibility. This paper extends the concept of reciprocity to linear time-and-space-invariant (LTSI) systems, a class of infinite-dimensional systems with spatio-temporal dynamics. It is suggested that, analogously to the LTI case, combining the internal properties of reciprocity and (impedance) passivity entails physical state-space realizations. This is of particular relevance for infinite-dimensional systems, where issues of unboundedness can be detrimental to the well-posedness of the system. The results are motivated and illustrated with a physical example.
	\end{abstract}
	
	\section{Introduction}
	
	The property of reciprocity refers to a type of symmetry in physical systems. The term originates in physics, where it appears in many forms: the Rayleigh–Carson reciprocity in acoustics; the Maxwell–Betti reciprocal work theorem in elasticity; Green’s reciprocity in electrostatics; Lorentz reciprocity in electromagnetics; and Onsager’s reciprocal relations in thermodynamics. In the field of systems and control theory, reciprocity was formally defined for linear time-invariant (LTI) systems in the seminal paper \cite{willems1972b}, see also \cite{willems1976,hughes2019, vanderschaft2025}. There, it was shown that, despite being an external property, reciprocity constrains the system’s internal structure. In particular, reciprocal systems admit a potential function and can thus be expressed as pseudo-gradient systems, even in the nonlinear case \cite{vanderschaft2024}. Passive and reciprocal systems can be realized using only LTI capacitors, inductors, resistors, and transformers. Furthermore, they admit a port-Hamiltonian formulation with a clear distinction between different types of energy storage elements. This explicit structure offers significant advantages for analysis and control design \cite{vanderschaft2017, pates2022}.
	
	Special subclasses of reciprocal systems include time-reversible and relaxation systems. The former is a combination of reciprocity and losslessness, while the latter is a combination of reciprocity and positivity. Relaxation systems, in particular, have received renewed attention in recent years \cite{pates2019, chaffey2023, sepulchre2024}. They are classically defined as being devoid of any oscillatory behaviour, which is expressed as having a completely monotone impulse response. As a subclass of reciprocal systems, they are defined as passive systems with a single type of energy storage element. They also admit a unique compatible storage function that can be determined from the external behaviour of the system.
	
	In this paper, we study the property of reciprocity for a class of infinite-dimensional systems that depend on both time and space, namely, linear time-and-space-invariant (LTSI) systems \cite{curtain2020, bamieh2002, arbelaiz2024}. The relaxation property for such systems has already been studied in \cite{donchev2025arxiv}, where it was shown that the characteristic features of relaxation can be seamlessly extended to the spatio-temporal domain. The focus of the present paper is slightly different. While we also aim to extend the property of reciprocity and its implications for LTSI systems, our focus is on the combination of reciprocity and (impedance) passivity, with the goal of deriving a well-posed and physically meaningful state-space realization. Unlike in finite-dimensional systems, issues of unboundedness and domain specification play a central role in ensuring stability and solvability in infinite-dimensional systems. We demonstrate that, for LTSI systems, these issues are resolved when the external physical properties of reciprocity and passivity are built into the internal structure of the system.
	
	The outline of this paper is as follows. In Section~\ref{sec:example}, we present the model of the Timoshenko beam to illustrate and motivate the contents of this paper. In Section~\ref{sec:notation_preliminaries}, we introduce notation and discuss relevant preliminaries. In Section~\ref{sec:ltsi_reciprocity}, we define reciprocity for LTSI systems and provide insight into the internal structure of reciprocal LTSI systems. In Section~\ref{sec:ltsi_reciprocity_passivity}, we combine reciprocity and impedance passivity for LTSI systems to arrive at the main result of this paper. Finally, we end the paper with concluding remarks in Section~\ref{sec:conclusion}.

	\section{Motivating example}\label{sec:example}
	
	In this section, we present a simple physical example to illustrate and motivate the contents of this paper. We consider the partial differential equations describing the dynamics of the Timoshenko beam \cite{graff1991}. These are given by
	\begin{equation}\label{eq:timoshenko_pde}
		\begin{aligned}
			\frac{\partial^2 q}{\partial t^2}(t,x) &= \frac{\partial^2 q}{ \partial x^2}(t,x)  - \frac{\partial \varphi}{ \partial x}(t,x) + u(t,x),\\
			\frac{\partial^2 \varphi}{\partial t^2}(t,x) &= - \frac{\partial q}{ \partial x}(t,x) + \frac{\partial^2 \varphi}{ \partial x^2}(t,x) - \varphi(t,x) ,
		\end{aligned}
	\end{equation}
	where $q(t,x)$ is the translational displacement, $\varphi(t,x)$ is the angular displacement, and $u(t,x)$ is the distributed load, each at time $t\in\bbR$ and position $x\in\bbR$. We assume an infinite homogeneous beam and we have set all constants to 1 for simplicity. Consequently, the dynamics of the Timoshenko beam are linear and invariant in both time and space.\phantom{.}

	The potential energy of the beam is the sum of the strain energies due to shear and bending, namely,
	\begin{equation*}
		P(t) = \half\int_{\bbR} \left(\frac{\partial q}{\partial x}(t,x) - \varphi(t,x)\right)^2 + \left(\frac{\partial \varphi}{\partial x} (t,x)\right)^2\, \d x.
	\end{equation*}
	The kinetic energy is
	\begin{equation*}
		K(t) = \half\int_{\bbR} \left(\frac{\partial q}{\partial t}(t,x)\right)^2 + \left(\frac{\partial \varphi}{\partial t}(t,x)\right)^2 \d x,
	\end{equation*}
	and the total energy is, thus, $E(t) = P(t) + K(t)$. Under suitable boundary conditions, it can be verified that
	\begin{equation*}
		\frac{\d E}{\d t}(t) = \int_{\bbR} \frac{\partial q}{\partial t}(t,x) u(t,x)\, \d x,
	\end{equation*}
	which suggests that the Timoshenko beam defines a passive (lossless) system with the distributed load $u(t,x)$ as input and the translational velocity $y(t,x) = \frac{\partial q}{\partial t}(t,x)$ as output. 
	
	This can also be recognized from a purely input-output perspective. To this end, take the spatial Fourier transform of \eqref{eq:timoshenko_pde} to obtain the family of ordinary differential equations
	\begin{equation}\label{eq:timoshenko_ode_omega}
		\begin{aligned}
			\frac{\d^2}{\d t^2}\hat q(t,\omega)&= - \omega^2 \hat q(t,\omega) + j\omega\hat \varphi(t,\omega) + \hat u(t,\omega),\\
			\frac{\d^2}{\d t^2}\hat \varphi(t,\omega)&= - j\omega \hat q(t,\omega) -\omega^2\hat \varphi(t,\omega) - \hat \varphi(t,\omega) ,
		\end{aligned}
	\end{equation}
	parameterized by the frequency $\omega\in\bbR$, where $\hat q(t,\cdot)$, $\hat \varphi(t,\cdot)$, and $\hat u(t,\cdot)$ are the Fourier transforms of $q(t,\cdot)$, $\varphi(t,\cdot)$, and $u(t,\cdot)$. Let $\hat y(t,\cdot)$ be the Fourier transform of $y(t,\cdot)$. Then, the temporal Laplace transform of \eqref{eq:timoshenko_ode_omega} provides the transfer function from $\hat u(\cdot,\omega)$ to $\hat y(\cdot,\omega)$:   
	\begin{equation}\label{eq:timoshenko_G_omega}
		\hat G(s,\omega) = \frac{s^3 + (\omega^2 +1)s}{s^4 + (2\omega^2+1)s^2 + \omega^4}.
	\end{equation}
	The latter is positive real for all $\omega\in\bbR$, hence \eqref{eq:timoshenko_ode_omega} defines a one-parameter family of passive systems, which again suggests that the Timoshenko beam defines a passive system.
	
	For finite-dimensional systems, a positive real transfer function guarantees the existence of a well-posed, passive state-space realization. For infinite-dimensional systems, however, the situation is more delicate. Although a system’s transfer function may be positive real for each frequency, naive choices of state variables can yield ill-posed dynamics, for example when the state operator is not the generator of a $C_0$-semigroup. Even if the system is well-posed, the associated storage functional need not be bounded on the state space.
	
	To illustrate this more precisely, note that, for every $\omega\in\bbR$, we can obtain a state-space realization of the finite-dimensional system \eqref{eq:timoshenko_ode_omega} by choosing the state as
	\begin{equation}\label{eq:timoshenko_state_omega}
		\hat z(t,\omega) = \bbm \hat q(t,\omega) & \hat \varphi(t,\omega) & \frac{\d \hat q}{\d t}(t,\omega) & \frac{\d \hat \varphi}{\d t}(t,\omega) \ebm\t.
	\end{equation}
	This yields a well-posed and passive LTI system, i.e., there exists a bounded positive definite storage function of the state that satisfies the associated dissipation inequality. The analogous state for the infinite-dimensional system \eqref{eq:timoshenko_pde} is 
	\begin{equation}
		z(t) = \bbm q(t,\cdot)& \varphi(t, \cdot) &  \frac{\partial q}{\partial t}(t, \cdot) & \frac{\partial \varphi}{\partial t}(t, \cdot) \ebm\t,
	\end{equation}
	Using this, we can write \eqref{eq:timoshenko_pde} as the LTSI system
	\begin{equation}\label{eq:timoshenko_ss_naive}
		\frac{\d}{\d t} z(t) = \calA z(t) + \calB u(t),
	\end{equation}
	where $u(t) = u(t,\cdot)$ and
	\begin{equation}\label{eq:timoshenko_AB}
		\calA = \bbm 0 & 0 & 1 & 0 \\
		0 & 0 & 0 & 1 \\
		\frac{\partial^2}{ \partial x^2} &  -\frac{\partial}{ \partial x} & 0 & 0 \\
		\frac{\partial}{ \partial x} & \frac{\partial^2}{ \partial x^2} - 1 & 0 & 0 \ebm, \quad \calB = \bbm 0 \\ 0 \\ 1 \\ 0 \ebm.
	\end{equation}
	
	It can be shown that the operator $\calA$ is \emph{not} the generator of a $C_0$-semigroup, hence \eqref{eq:timoshenko_ss_naive} is not well-posed. Furthermore, the total energy is $E(t) = \half\inner{\calQ z(t)}{z(t)}$, where $\inner{\cdot}{\cdot}$ is the inner product in $L_2(\bbR,\bbR^4)$ and
	\begin{equation}\label{eq:calQ}
		\calQ = \bbm -\frac{\partial^2}{\partial x^2} & \frac{\partial}{\partial x} & 0 & 0 \\  -\frac{\partial}{\partial x} &  1 - \frac{\partial^2}{\partial x^2} & 0 & 0 \\ 0 & 0 & 1 & 0 \\ 0 & 0 & 0 & 1 \ebm.
	\end{equation}
	Clearly, the operator $\calQ$ that defines the total energy is not bounded. Although one can define and study passivity with unbounded storage \cite{reis2025}, these issues suggest that the choice of state is not appropriate. The total energy rather suggests that the appropriate physical states are the shear force $\frac{\partial q}{\partial x}(t,\cdot) - \varphi(t,\cdot)$ and the bending moment $\frac{\partial \varphi}{\partial x}(t, \cdot)$ rather than the translational displacement $q(t,\cdot)$ and the angular displacement $\varphi(t,\cdot)$. In fact, the Timoshenko beam can be formulated as a port-Hamiltonian system  with state
	\begin{equation}\label{eq:timoshenko_barz}
		\bar z(t) = \big[ \tfrac{\partial q}{\partial x}(t,\cdot) - \varphi(t,\cdot)\ \ \tfrac{\partial \varphi}{\partial x}(t, \cdot)\ \ \tfrac{\partial q}{\partial t}(t, \cdot)\ \  \tfrac{\partial \varphi}{\partial t}(t, \cdot) \big]\t,
	\end{equation}
	see \cite[Example~7.1.4]{jacob2012} for details.
	
	With this in mind, the main question we address in this paper is the following: how can one construct a well-posed and physically meaningful state-space realization of an LTSI system? We demonstrate that incorporating the external properties of reciprocity and passivity into the internal structure of the system is essential. The key idea in our analysis is that LTSI systems can be decomposed into an infinite family of LTI systems via the Fourier transform, an idea that is at the core of the LTSI property and has already been exploited in, e.g., \cite{bamieh2002, arbelaiz2024, donchev2025arxiv}.

	\section{Notation and preliminaries}\label{sec:notation_preliminaries}
	
	We use mostly standard notation. We denote the set of nonnegative real numbers by $\bbRp$. We denote the inner product on a Hilbert space $\calX$ by $\inner{\cdot}{\cdot}_\calX$ and the induced norm by $\norm{\cdot}_\calX$. We omit the subscript if there is no risk of ambiguity. We denote the Banach space of bounded linear operators from $\calX$ to another Hilbert space $\calY$ by $\B(\calX,\calY)$ and the operator norm by $\norm{\cdot}$. We use the shorthand notation $\B(\calX,\calX) = \B(\calX)$. We also consider unbounded linear operators $\calA:\calX\supset\D(\calA) \to \calY$ with domain $\D(\calA)$. We say that $\calA$ is densely defined if $\D(\calA)$ is dense. Every densely defined $\calA$ has an adjoint $\calA\a: \calY\supset\D(\calA\a) \to \calX$ that is a \emph{closed} operator \cite[Appendix~A.3]{curtain2020}. We say that $\calA$ is self-adjoint if $\D(\calA) = \D(\calA\a)$ and $\inner{\calA x}{y} = \inner{x}{\calA y}$ for all $x,y\in\D(\calA)$. We say that a self-adjoint $\calA$ is nonnegative if $\inner{\calA x}{x} \geq 0$ for all $x\in\D(\calA)$, and we write $\calA \geq 0$ to indicate that $\calA$ is self-adjoint and nonnegative. We use similar notation to indicate that an operator is self-adjoint and nonpositive, positive, or negative.
	
	\subsection{Square integrable functions and the Fourier transform}
	We denote the Hilbert space of square integrable functions from a domain $X$ to a Hilbert space $\calX$ by $L_2(X, \calX)$, where
	\begin{equation}
		\inner{f}{g}_{L_2(X,\calX)} = \int_{X} \inner{f(x)}{g(x)}_\calX\, \d x.
	\end{equation}
	The Fourier transform of a function $f:\bbR^s\to\bbC^n$ is the function $\hat f: \bbR^s\to\bbC^n$ given by
	\begin{equation}
		\hat f(\omega) = (2\pi)^{-s/2}\intrs e^{-j\inner{\omega}{x}}f(x)\, \d x.
	\end{equation}
	It is well-known \cite{zuazo2001} that the Fourier transform converges for absolutely integrable functions. By Plancherel's theorem, the Fourier transform can be extended to a \emph{unitary} operator $\calF\in \B(L_2(\bbR^s, \bbC^n))$, i.e., $\calF\inv = \calF\a$. Since $\calF$ agrees with the Fourier transform on the dense subspace of absolutely integrable functions, we refer to $\calF f$ as the Fourier transform of $f$, and we write $\hat f = \calF f$. 
	
	Suppose that $\calA:L_2(\bbR^s,\bbC^n) \supset \D(\calA)\to L_2(\bbR^s,\bbC^m)$ is densely defined. The Fourier transform of $\calA$ is the densely defined operator $\hat \calA: L_2(\bbR^s,\bbC^n) \supset \D(\hat \calA)\to L_2(\bbR^s,\bbC^m)$ such that $\hat \calA \hat f$ is the Fourier transform of $\calA f$. For $x\in\bbR^s$, the \emph{translation operator} $\calT_{x}\in\B(L_2(\bbR^s,\bbC^n))$ is given by $\calT_{x} f = f(\cdot+ x)$. We say that $\calA$ is \emph{translation invariant} if $\calT_x(\D(\calA)) \subset \D(\calA)
	$ and $\calT_{x}\calA = \calA \calT_{x}$ for all $x\in\bbR^s$. The Fourier transform of a translation invariant operator is a \emph{multiplication operator} \cite[Section~3.6.7]{zuazo2001}, i.e., there exists a function $\hat A:\bbR^s\to \bbC^{m\times n}$ such that $(\hat\calA \hat f)(\omega) = \hat A(\omega) \hat f(\omega)$. The function $\hat A$ is called the \emph{symbol} of the multiplication operator $\hat\calA$. In this paper, we also refer to $\hat A$ as the \emph{Fourier symbol} of the translation invariant operator $\calA$. It is well-known that $\hat \calA$ and, thus, $\calA$ are bounded if and only if $\hat A$ is bounded. Conversely, every function $\hat A:\bbR^s\to\bbC^{m\times n}$ that is bounded on compact subsets of $\bbR^s$ (e.g., any continuous $\hat A$) yields a densely defined translation invariant operator $\calA$ with Fourier symbol $\hat A$. It is easily seen that $\calA\a$ is the translation invariant operator with Fourier symbol $\omega\mapsto \hat A(\omega)\a$, hence $\calA\a$ is also densely defined.

	Throughout this paper, we consider functions of the form $f:\bbT\to  L_2(\bbR^s,\bbC^n)$, where $\bbT$ is a temporal domain and $\bbR^s$ is a spatial (frequency) domain. We denote the dependence on the spatial (frequency) variable with a subscript, i.e., we write $f_x(t) = f(t)(x)$, so that the function $f_x:\bbT\to \bbC^n$ is the temporal trajectory at a given point in space $x\in\bbR^s$. Finally, $\hat f$ denotes the spatial Fourier transform of $f$, i.e., the function $\hat f: \bbT\to  L_2(\bbR^s,\bbC^n)$ such that $\hat f(t)$ is the Fourier transform of $f(t)$ for all $t\in\bbT$.
	
	\subsection{Reciprocal LTI systems}
	
	Consider the linear time-invariant (LTI) system
	\begin{equation}\label{eq:lti_ss}
		\sys:\left\lbrace \begin{aligned}
			\frac{\d }{\d t}z(t) &= Az(t) + Bu(t),\\
			y(t) &= Cz(t),
		\end{aligned}\right.
	\end{equation}
	with input space $\bbR^m$, state space $\bbR^n$, output space $\bbR^p$, and matrices $A\in\bbR^{n\times n}$, $B\in\bbR^{n\times m}$, and $C\in\bbR^{p\times n}$. The \emph{impulse response} of $\sys$ is the map $g:\bbRp\to\bbR^{p\times m}$ given by $g(t) = Ce^{At} B$. The \emph{transfer matrix} of $\sys$ is the map $G:\bbC\to \bbC^{p\times m}$ given by $G(s) = C(sI-A)\inv B$. The output of $\sys$ is given by the convolution integral
	\begin{equation}
		y(t) = \int_{-\infty}^{t} g(t-\tau)u(\tau)\, \d \tau,
	\end{equation}
	where $z(-\infty) = 0$, i.e., the system is initially at rest.

	We are interested in systems $\sys$ that are \emph{reciprocal} \cite{willems1972b}. Intuitively, $\sys$ is reciprocal if $m=p$ and, for all $i,j\in\{1,\dots,m\}$, the effect of the $i$'th input on the $j$'th output is the same as the effect of the $j$'th input on the $i$'th output, see, e.g., \cite[Definition~4]{willems1976}. In general, reciprocity is defined with respect to a signature matrix. For simplicity, we only consider reciprocity with respect to the identity matrix and we define it in terms of its characterization, namely, $\sys$ is reciprocal if $g(t) = g(t)\t$ for all $t\in\bbRp$, or, equivalently, $G(s) = G(s)\t$ for all $s\in\bbC$. 
	
	Reciprocity can also be characterized internally. In particular, a minimal $\sys$ is reciprocal if and only if there exists a (unique) invertible symmetric matrix $S\in\bbR^{n\times n}$ such that
	\begin{equation}\label{eq:lti_S}
		A\t S = SA ,\qquad C\t = SB.
	\end{equation}
	This implies that $\sys$ is \emph{self-dual} \cite{vanderschaft2025}, i.e., it is isomorphic to the dual LTI system
	\begin{equation}\label{eq:lti_ss_dual}
		\sys\dual:\left\lbrace \begin{aligned}
			\frac{\d }{\d t}z\dual(t) &= A\t z\dual(t) + C\t u\dual(t),\\
			y\dual(t) &= B\t z\dual(t),
		\end{aligned}\right.
	\end{equation}
	with $z\dual = Sz$ being the unique state-space isomorphism. If $\sys$ is asymptotically stable, then $S$ defines a potential that can be determined from the input-output behaviour, namely,
	\begin{equation}\label{eq:lti_potential_external}
		\half z(0)\t Sz(0) = \int_{0}^{\infty}y(t)u(-t)\, \d t.
	\end{equation}
	The quantity on the left-hand side is called the \emph{Lagrangian} and is invariant under state-space transformations. There exists a state-space transformation in which $\sys$ is internally  reciprocal with respect to the signature of $S$, i.e., 
	\begin{equation}\label{eq:partition}
		A = \bbm A_{11} & A_{12} \\ -A_{12}\t & A_{22} \ebm, \quad B = \bbm B_1 \\ B_2 \ebm,\quad C\t = \bbm B_1 \\ - B_2 \ebm 
	\end{equation}
	with $A_{11} = A_{11}\t$ and $A_{22} = A_{22}\t$. 
	
	Recall that a minimal $\sys$ is \emph{passive} if and only if there exists $Q\in\bbR^{n\times n}$, $Q>0$, such that
	\begin{equation}
		A\t Q + QA \leq 0,\quad C\t = QB
	\end{equation}
	with $\half x\t Q x$ being the corresponding storage function \cite{willems1972b}. We say that $\sys$ is \emph{internally passive} if it is passive with storage $\half\norm{x}^2$, i.e., $A\t + A \leq 0$ and $C\t = B$. If $\sys$ is both reciprocal and passive, then there exists a so-called \emph{compatible} storage, i.e., $Q>0$ such that $Q = SQ\inv S$. Consequently, there exists a state-space transformation that is both internally passive and reciprocal with respect to the signature of $S$, i.e., \eqref{eq:partition} holds with $A_{11} \leq 0$, $A_{22}\leq 0$ and $B_2 = 0$. The latter has the structure of a linear port-Hamiltonian system. Partitioning the state according to \eqref{eq:partition} shows that a passive reciprocal system contains two types of energy storage elements, e.g., capacitors and inductors, between which there is lossless energy exchange.
	
	If $S$ is positive definite, then $\sys$ is of \emph{relaxation} type \cite{willems1972b}, i.e., its impulse response is completely monotone \cite{bernstein1929}. If, in addition, $\sys$ is passive, then $Q = S$ defines the unique compatible storage, which, due to \eqref{eq:lti_potential_external}, can be determined from input-output behaviour. Moreover, there exists a state-space transformation in which $\sys$ is both internally passive and \emph{symmetric}, i.e., $A \leq 0$ and $C\t = B$. This shows that relaxation systems contain only one type of energy storage element, e.g., only capacitors.
	
	The above definition of reciprocity is for real systems, i.e., $\sys$ of the form \eqref{eq:lti_ss} with real input, state, and output spaces, as well as real system matrices $A$, $B$, and $C$. However, throughout this paper, we have to consider complex systems, which have, in particular, complex impulse responses. With this in mind, a complex $\sys$ is reciprocal if $g(t) = g(t)\a$ for all $t\in\bbRp$, where $g(t)\a\in\bbC^{m\times m}$ is the adjoint (conjugate transpose) of $g(t)\in\bbC^{m\times m}$. The latter is equivalent to $G(s)\a = G(s\a)$ for all $s\in\bbC$. The internal properties of real reciprocal systems are also valid for complex reciprocal systems. In particular, if $\sys$ is minimal, there exists an invertible self-adjoint matrix $S\in\bbC^{n\times n}$ such that \eqref{eq:lti_S} holds with transposes replaced by adjoints.  If $\sys$ is also asymptotically stable, then \eqref{eq:lti_potential_external} also holds. If $\sys$ is also passive, then there exists a compatible storage and, thus, a state-space transformation in which \eqref{eq:partition} holds with transposes replaced by adjoints, $A_{11}\leq 0$, $A_{22}\leq 0$ and $B_2 = 0$.

	\section{Reciprocal LTSI systems}\label{sec:ltsi_reciprocity}
	
	In this section, we define and characterize the notion of reciprocity for a class of linear time-and-space-invariant (LTSI) systems. The definition is a natural generalization of the definition for LTI systems, that is, it is expressed as self-adjointness of the (operator-valued) impulse response. Due to the spatial invariance of the LTSI system, we can decouple it into a family of LTI systems via the Fourier transform. Consequently, we can characterize reciprocity of the LTSI system in terms of reciprocity of the associated family of LTI systems. This characterization allows us to extend some of the internal properties of reciprocal LTI systems to internal properties of reciprocal LTSI systems.
	
	Consider the LTSI system $\sys$ of the form
	\begin{equation}\label{eq:ltsi_ss}
		\sys:\left\lbrace\begin{aligned}
			\frac{\d}{\d t}z(t) &= \calA z(t) + \calB u(t),\\
			y(t) &= \calC z(t),
		\end{aligned}\right.
	\end{equation}
	with (infinite-dimensional) input space $\calU = L_2(\bbR^s, \bbC^m)$, state space $\calZ = L_2(\bbR^s, \bbC^n)$, output space $\calY = L_2(\bbR^s, \bbC^p)$, and translation invariant operators $\calA:\calZ\supset\calD(\calA) \to \calZ$, $\calB \in \B(\calU,\calZ)$ and $\calC\in\B(\calZ,\calY)$, where we assume that the domain $\D(\calA)$ is dense. The system is time invariant because the operators $\calA$, $\calB$, and $\calC$ are constant, and space-invariant because they are translation invariant. The operator $\calA$ is typically unbounded and the notion of solution requires some care, see, e.g. \cite{curtain2020}. Without going into details, we say that $\sys$ is well-posed if $\calA$ is the infinitesimal generator of a $C_0$-semigroup $t\mapsto \exp(\calA t)\in\B(\calZ)$, where the notation is motivated by the fact that the latter is a generalization of the exponential function.
	
	The impulse response of a well-posed LTSI system $\sys$ is the operator-valued map $\calG:\bbRp\to\B(\calU,\calY)$ given by
	\begin{equation}
		\calG(t) = \calC\exp(\calA t)\calB.
	\end{equation}
	The output of $\sys$ is given by the convolution
	\begin{equation}
		y(t) = \int_{-\infty}^{t} \calG(t-\tau)u(\tau)\, \d\tau,
	\end{equation}
	where we have assumed that $z(-\infty) = 0$, i.e., the system is initially at rest. Since $\calA$, $\calB$ and $\calC$ are translation invariant, their Fourier transforms $\hat \calA$, $\hat\calB$, and $\hat\calC$ are multiplication operators with symbols $\hat A$, $\hat B$ and $\hat C$. Consequently, the Fourier transform of \eqref{eq:ltsi_ss} yields the family of LTI systems
	\begin{equation}\label{eq:lti_ss_omega}
		\hat\sys_\omega: \left\lbrace
		\begin{aligned}
			\frac{\d}{\d t}\hat z_\omega(t) &= \hat A_\omega\hat z_\omega (t) + \hat B_\omega \hat u_\omega(t), \\ \hat y_\omega(t)&= \hat C_\omega\hat z_\omega(t),
		\end{aligned}\right.
	\end{equation}
	parametrized by the frequency variable $\omega\in\bbR^s$. The Fourier transform $\hat \calG(t)$ of $\calG(t)$ is a multiplication operator with symbol $\hat g(t)$ given by $	\hat g_\omega(t) = \hat C_\omega e^{\hat A_\omega t}\hat B_\omega$ for all $\omega\in\bbR^s$, where $\hat g_\omega$ is the impulse response of the LTI system $\hat\sys_\omega$. We assume that the symbols $\hat A$, $\hat B$ and $\hat C$ are continuous, so that $\hat g(t)$ is continuous for all $t\in\bbRp$.
	
	Even if $\sys$ is not well-posed, the Fourier transform of \eqref{eq:ltsi_ss} still yields a family of LTI systems \eqref{eq:lti_ss_omega}. Consequently, for all $t\in\bbRp$, we can define $\calG(t):\calU\supset\D(\calG(t))\to\calY$ as the translation invariant operator with Fourier symbol $\hat g(t)$. Since $\hat g(t)$ is continuous, $\D(\calG(t))$ is dense, but we are not guaranteed that $\D(\calG(t)) = \calU$ or that $\calG(t)$ is bounded. Nevertheless, we still treat $\calG$ as the impulse response of $\sys$ and use it for the definition of reciprocity.
	
	\begin{definition}
		The LTSI system $\sys$ is reciprocal if $\calU = \calY$ and $\calG(t)$ is self-adjoint for all $t\in\bbRp$.
	\end{definition}
	
	We can characterize reciprocity of an LTSI system as reciprocity of the associated family of LTI systems.
	
	\begin{theorem}\label{thm:LTSI_reciprocity_LTI}
		The LTSI system $\sys$ is reciprocal if and only if the LTI system $\hat \sys_\omega$ is reciprocal for all $\omega\in\bbR^s$.
	\end{theorem}
	\proof
	Since the Fourier transform is unitary, $\calG(t)$ is self-adjoint if and only if $\hat\calG(t)$ is self-adjoint. Since $\hat\calG(t)$ is a multiplication operator with continuous symbol $\hat g(t)$, it is self-adjoint if and only if $\hat g_\omega(t)$ is self-adjoint for all $\omega\in\bbR^s$. Therefore, the LTSI system $\sys$ is reciprocal if and only if the LTI system $\hat \sys_\omega$ is reciprocal for all $\omega\in\bbR^s$.
	\endproof
	
	We illustrate the definition and characterization of reciprocity with the example of the Timoshenko beam.
	\begin{example}
		As shown in Section~\ref{sec:example}, the Timoshenko beam can be formulated as an (ill-posed) LTSI system $\sys$ of the form \eqref{eq:ltsi_ss} with $\calA$ and $\calB$ as in \eqref{eq:timoshenko_AB}, and $\calC = \calB\a$. Using the properties of the Fourier transform, $\sys$ is transformed into a family of LTI systems $\hat \sys_\omega$, $\omega\in\bbR$, of the form \eqref{eq:lti_ss_omega} with
		\begin{equation*}
			\hat A_\omega = 
			\bbm 0 & 0 & 1 & 0 \\
			0 & 0 & 0 & 1 \\
			-\omega^2 & j\omega & 0 & 0 \\
			-j\omega & -\omega^2 - 1 & 0 & 0 \ebm,
			\quad
			\hat C_\omega\a = \hat B_\omega = \bbm 0 \\ 0 \\ 1 \\ 0 \ebm.
		\end{equation*}
		The transfer function of $\hat \sys_\omega$ is given by \eqref{eq:timoshenko_G_omega}, which clearly satisfies $\hat G_\omega(s)\a = \hat G_\omega(s\a)$ for all $s\in\bbC$. Therefore, $\hat \sys_\omega$ is reciprocal for all $\omega\in\bbR$ and, due to Theorem~\ref{thm:LTSI_reciprocity_LTI}, the LTSI system $\sys$ describing the Timoshenko beam is reciprocal.
	\end{example}
	
	In the following, we provide an internal characterization of reciprocity for the LTSI system $\sys$. We do this under the assumption that the LTI system $\hat \sys_\omega$ is minimal for all $\omega\in\Omega$, where $\Omega\subset\bbR^s$ is dense. If $\sys$ is well-posed, this is equivalent to $\sys$ being approximately controllable and observable (minimal) on any time interval, see \cite[Definition~6.2.1]{curtain2020} and \cite[Theorem~6.3.1]{curtain2020}. Under this assumption, the controllability and observability matrices of $\hat \sys_\omega$, given by
	\begin{align}
		\hat W_\omega &= \bbm \hat B_\omega & \hat A_\omega \hat B_\omega & \cdots & \hat A_\omega^{n-1} \hat B_\omega \ebm, \\
		\hat O_\omega\a &= \bbm \hat C_\omega\a & \hat A_\omega\a \hat C_\omega\a & \cdots & (\hat A_\omega\a)^{n-1} \hat C_\omega\a \ebm,
	\end{align}
	respectively, have full rank for all $\omega\in\Omega$. Consequently, we obtain the following theorem.	
	
	\begin{theorem}\label{thm:calS}
		Suppose that the LTSI system $\sys$ is reciprocal, the LTI system $\hat \sys_\omega$ is minimal for all $\omega\in\Omega$, where $\Omega\subset\bbR^s$ is dense, and
		\begin{equation}\label{eq:S_omega_continuiuty}
			\lim_{\omega\to \bar\omega, \omega\in\Omega} \hat O_\omega\a \hat W_\omega\a (\hat W_\omega \hat W_\omega\a)\inv
		\end{equation}
		exists for all $\bar\omega\in\bbR^s\setminus\Omega$. Then, there exists a unique, densely defined, self-adjoint, translation invariant linear operator\linebreak $\calS:\calZ\supset\D(\calS) \to\calZ$, whose Fourier symbol $\hat S:\bbR^s\mapsto\bbC^{n\times n}$ is continuous and satisfies
		\begin{equation}\label{eq:S_omega_definition}
			\hat S_\omega = \hat O_\omega\a \hat W_\omega\a (\hat W_\omega \hat W_\omega\a)\inv
		\end{equation}
		for all $\omega\in\Omega$, and
		\begin{equation}\label{eq:calS}
			\inner{\calA z_1}{\calS z_2}= \inner{\calS z_1}{\calA z_2},\quad \calC\a = \calS \calB.
		\end{equation}
		for all $z_1,z_2\in\D(\calA)\cap\D(\calS)$ and $u\in\calU$.
	\end{theorem}
	
	\proof
		By Theorem~\ref{thm:LTSI_reciprocity_LTI}, the assumption that $\sys$ is reciprocal implies that $\hat \sys_\omega$ is reciprocal for all $\omega\in\bbR^s$. Since $\hat \sys_\omega$ is minimal and reciprocal for all $\omega\in\Omega$, there exists a unique invertible self-adjoint matrix $\hat S_\omega\in\bbC^{n\times n}$ such that 
		\begin{equation}\label{eq:S_omega_property}
			\hat A_\omega\a \hat S_\omega = \hat S_\omega \hat A_\omega, \quad \hat C_\omega \a = \hat S_\omega \hat B_\omega
		\end{equation}
		for all $\omega\in\Omega$. Note that $\hat S_\omega \hat W_\omega = \hat O_\omega\a$, hence \eqref{eq:S_omega_definition} holds for all $\omega\in\Omega$. Recall that the symbols  $\hat A$, $\hat B$ and $\hat C$ are continuous, hence $\omega\mapsto (\hat W_\omega, \hat O_\omega)$ is continuous on $\bbR^s$. This implies that the map $\omega\mapsto \hat S_\omega$ is continuous on $\Omega$, where the inverse of $\hat W_\omega \hat W_\omega\a$ exists. Due to \eqref{eq:S_omega_continuiuty}, we can extend this to a map $\omega\mapsto \hat S_\omega$ that is continuous on the whole of $\bbR^s$. Since both sides of the equations in \eqref{eq:S_omega_property} are continuous in $\omega$, \eqref{eq:S_omega_property} holds for all $\omega\in\bbR^s$. Consequently, the translation invariant operator $\calS:\calZ\supset\D(\calS) \to\calZ$ with Fourier symbol $\omega\mapsto \hat S_\omega$ is densely-defined and self-adjoint, while \eqref{eq:S_omega_property} implies that $\im\calB\subset \D(\calS)$ and \eqref{eq:calS} holds.
	\endproof
	
	We illustrate the use of Theorem~\ref{thm:calS} with the example of the Timoshenko beam.
	\begin{example}\label{ex:timoshenko_reciprocal}
		For the (ill-posed) LTSI system $\sys$ that describes the Timoshenko beam, we have that 
		\begin{equation}
			\hat W_\omega = \bbm 0 & 1& 0 & -\omega^2 \\
			0 & 0 & 0 & -j\omega \\
			1 & 0 & -\omega^2 & 0 \\
			0 & 0 & -j\omega & 0 \ebm,
		\end{equation}
		hence $\hat \sys_\omega$ is not controllable for $\omega = 0$. Similarly, we can show that $\hat \sys_\omega$ is not observable for $\omega = 0$. However, $\hat \sys_\omega$ is both controllable and observable for all $\omega \neq 0$, which is a dense subset of $\bbR$. Furthermore, \eqref{eq:S_omega_definition} reduces to
		\begin{equation}\label{eq:timoshenko_S_omega}
			\hat S_\omega =
			\bbm -\omega^2 & j\omega & 0 & 0 \\ 
			-j\omega & -\omega^2 - 1 & 0 & 0 \\ 
			0 & 0 & 1 & 0 \\ 
			0 & 0 & 0 & 1 \ebm
		\end{equation}
		where we observe that $\omega\mapsto\hat S_\omega$ is continuous on $\bbR$. It is straightforward to verify that \eqref{eq:S_omega_property} holds for all $\omega\in\bbR$ and the translation invariant operator $\calS$ with Fourier symbol $\omega\mapsto\hat S_\omega$ satisfies $\im\calB\subset \D(\calS) = \D(\calA)$ and \eqref{eq:calS}.
	\end{example}
	
	Recall that LTI systems that are minimal and reciprocal are also self-dual. This is not necessarily the case for LTSI systems, mainly due to issues of unboundedness. To see this, note that the dual of the LTSI system $\sys$ is the LTSI system
	\begin{equation}\label{eq:ltsi_ss_d}
		\sys\dual:\left\lbrace\begin{aligned}
			\frac{\d}{\d t}z^\d(t) &= \calA\a z^\d(t) + \calC\a u\dual(t),\\
			y\dual(t) &= \calB\a z\dual(t),
		\end{aligned}\right.
	\end{equation}
	with input space $\calU\dual = \calY$, state space $\calZ\dual = \calZ$ and output space $\calY\dual = \calU$. The spatial Fourier transform of \eqref{eq:ltsi_ss_d} yields the family of LTI systems $\hat\sys_\omega\dual$, $\omega\in\bbR^s$, where $\hat\sys_\omega\dual$ is the dual of $\hat\sys_\omega$ for all $\omega\in\bbR^s$. If $\calA$ is the infinitesimal generator of a $C_0$-semigroup, then $\calA\a$ is as well \cite[Theorem~2.3.6]{curtain2020}. This implies that $\sys\dual$ is well-posed if and only if $\sys$ is well-posed, and its impulse response $\calG\dual:\bbRp\to \B(\calY,\calU)$ is given by
	\begin{equation}
		\calG\dual(t) = \calB\a\exp(\calA\a t)\calC\a = \calG(t)\a,
	\end{equation}
	where we used the fact that $\exp(\calA\a t) = \exp(\calA t)\a$. 
	
	Now, suppose that $\sys$ is well-posed and reciprocal. This implies that $\sys$ and $\sys\dual$ have the same impulse response, i.e., they are externally equivalent. Suppose further that $\sys$ is exactly minimal, i.e., $\hat \sys_\omega$ is minimal for all $\omega\in\bbR^s$. Due to Theorem~\ref{thm:LTSI_reciprocity_LTI}, the LTI system $\hat \sys_\omega$ is minimal and reciprocal, hence self-dual, that is, isomorphic to $\sys_\omega\dual$ for all $\omega\in\bbR^s$. This does not necessarily imply that $\sys$ is isomorphic to $\sys\dual$. Indeed, the unique state-space isomorphism that relates $\hat \sys_\omega$ to $\hat \sys_\omega\dual$ is $\hat z_\omega\dual = \hat S_\omega \hat z_\omega$, where $\hat S_\omega$ is given by \eqref{eq:S_omega_definition}. Therefore, the potential state-space isomorphism that relates $\sys$ to $\sys\dual$ must be given by $z\dual = \calS z$, where $\calS$ is the translation invariant operator with Fourier symbol $\omega\mapsto\hat S_\omega$. However, for this to constitute a genuine state-space isomorphism, we need $\calS$ to be bounded with bounded inverse, neither of which are guaranteed. Nevertheless, $\hat S_\omega$ is invertible and 
	\begin{equation}
		\hat S_\omega\inv = (\hat O_\omega \hat O_\omega\a)\inv\hat O_\omega\a \hat W_\omega\a
	\end{equation}
	for all $\omega\in\bbR^s$. With this in mind, we conclude this section with the following corollary of Theorem~\ref{thm:calS}.
	\begin{corollary}\label{cor:selfdual}
		Suppose that the well-posed LTSI system $\sys$ is reciprocal, and exactly controllable and observable. If
		\begin{align}
			\sup_{\omega\in\bbR^s} \norm{\hat O_\omega\a \hat W_\omega\a (\hat W_\omega \hat W_\omega\a)\inv} &< \infty,\label{eq:S_omega_bounded}\\
			\sup_{\omega\in\bbR^s} \norm{(\hat O_\omega \hat O_\omega\a)\inv\hat O_\omega\a \hat W_\omega\a} &< \infty\label{eq:S_omega_inv_bounded},
		\end{align}
		then $\sys$ is self-dual, i.e., it is isomorphic to $\sys\dual$ with $z\dual = \calS z$ being the unique state-space isomorphism, where $\calS\in\B(\calZ)$ is the self-adjoint translation invariant linear operator with Fourier symbol $\hat S:\bbR^s\mapsto\bbC^{n\times n}$ given by \eqref{eq:S_omega_definition}.
	\end{corollary}
	\proof
	 Note that \eqref{eq:S_omega_bounded} and \eqref{eq:S_omega_inv_bounded} ensure that $\calS$ is bounded with bounded inverse, hence $z\dual = \calS z$ is indeed a state-space isomorphism, as discussed above.
	\endproof
	
	\section{Impedance passive and reciprocal LTSI systems}\label{sec:ltsi_reciprocity_passivity}
	
	In this section, we combine the properties of reciprocity and impedance passivity for LTSI systems. Impedance passivity is the analogue of passivity for infinite-dimensional systems, defined via the existence of a suitable quadratic storage functional \cite{curtain2020}, see also \cite{staffans2002, arov2005}. This definition is for well-posed systems and requires the storage functional to be bounded, neither of which hold for the naive state-space realization of the Timoshenko beam in Section~\ref{sec:example}. The requirement on boundedness can be dropped \cite{reis2025}, but we still face the issue of well-posedness. To address that, we exploit the fact that an LTSI system can be seen as a family of LTI systems, so that a weak form of impedance passivity for the former can be defined by passivity of the latter. Then, we show that, under additional technical assumptions, reciprocity and this weak form of impedance passivity guarantee the existence of an externally equivalent LTSI system that is well-posed, (internally) impedance passive and self-dual.
	
	We begin with a definition of impedance passivity for well-posed LTSI systems, see \cite[Section~7.5]{curtain2020} for details.
	\begin{definition}
		The well-posed LTSI system $\sys$ is \emph{impedance passive} if $\calU = \calY$ and there exists a (quadratic) storage functional $\calS:\calZ\to\bbRp$, given by $\calS(z) = \inner{\calQ z}{z}$ with $\calQ\in\B(\calZ)$, $\calQ\geq 0$, such that 
		\begin{equation*}
			\calS(z(t)) \leq \calS(z(0)) + \int_{0}^{t}  \inner{u(\tau)}{y(\tau)} + \inner{y(\tau)}{u(\tau)}\, \d \tau
		\end{equation*}
		for all $z(0)\in\calZ$, $t\geq 0$, and $u\in L_2([0,\tau], \calU)$. We say that $\sys$ is \emph{internally impedance passive} if $\calS(z) = \half\norm{z}^2$.
	\end{definition}
	
	Impedance passivity can be characterized as the solvability of a linear operator inequality \cite[Lemma~7.5.4]{curtain2020}. For well-posed LTSI systems, it can also be characterized as passivity of the associated family of LTI systems. 
	
	\begin{proposition}[{\cite[Theorem~3]{donchev2025arxiv}}]\label{prop:impedance_passive}
		The well-posed LTSI system $\sys$ is impedance passive if and only if there exist $\hat Q_\omega\in\bbC^{n\times n}$, $\hat Q_\omega \geq 0$, such that
		\begin{equation}\label{eq:Q_omega_property}
			\hat A_\omega\a \hat Q_\omega + \hat Q_\omega \hat A_\omega \leq 0, \quad \hat C_\omega\a = \hat Q_\omega \hat B_\omega , 
		\end{equation} 
		for all $\omega\in\bbR^s$, and $\sup_{\omega\in\bbR^s} \norm{\hat Q_\omega} < \infty$. Furthermore, $\sys$ is internally impedance passive if and only if $\eqref{eq:Q_omega_property}$ holds with $\hat Q_\omega = I$ for all $\omega\in\bbR^s$.
	\end{proposition}
	
	As shown in Section~\ref{sec:example}, a naive state-space realization of an infinite-dimensional physical system can lead to an LTSI system that is neither well-posed nor impedance passive, although its external (physical) properties suggest otherwise. This motivates the following definition.
	\begin{definition}
		The LTSI system $\sys$ is \emph{weakly impedance passive} if the LTI system $\hat\sys_\omega$ is passive for all $\omega\in\bbR^s$.
	\end{definition}
	
	In other words, $\sys$ is weakly impedance passive if and only if there exists $\hat Q_\omega\in\bbC^{n\times n}$, $\hat Q_\omega \geq 0$, such that \eqref{eq:Q_omega_property} holds for all $\omega\in\bbR^s$. By Proposition~\ref{prop:impedance_passive}, a well-posed and impedance passive $\sys$ is automatically weakly impedance passive. The converse does not generally hold because the boundedness condition might fail, e.g., due to an inappropriate choice of state-space compromising well-posedness.
	In the following, we show that the external properties of reciproctity and weak impedance passivity enable the construction of an externally equivalent state-space realization that is well-posed, internally impedance passive and self-dual.
	
	To this end, suppose that $\sys$ is reciprocal and weakly impedance passive, so that $\hat \sys_\omega$ is reciprocal and passive for all $\omega\in\bbR^s$. We also assume that $\hat\sys_\omega$ is minimal for all $\bbR^s$. This is a stronger assumption than the one made in the previous section, where we required minimality of $\hat \sys_\omega$ only for $\omega$ in a dense subset of $\bbR^s$. The reason for this stronger assumption is mainly technical and we demonstrate that the results are useful even when it is not satisfied. 
	
	Now, the assumptions that $\hat\sys_\omega$ is reciprocal, passive, and minimal for all $\omega\in\bbR^s$ imply that the following hold for all $\omega\in\bbR^s$. First, there exists a unique invertible self-adjoint $\hat S_\omega\in\bbC^{n\times n}$ such that \eqref{eq:S_omega_property} holds. Note that $\hat S_\omega$ is given by \eqref{eq:S_omega_definition} and the map $\omega\mapsto\hat S_\omega$ is continuous on $\bbR^s$ because we have assumed that the symbols $\hat A$, $\hat B$ and $\hat C$ are continuous. Second, there exists a compatible storage function for $\hat \sys_\omega$, i.e., $\hat Q_\omega\in\bbC^{n\times n}$, $\hat Q_\omega > 0$, such $\hat Q_\omega = \hat S_\omega \hat Q_\omega\inv \hat S_\omega$ and \eqref{eq:Q_omega_property} holds. As shown in \cite[Theorem~7]{willems1972b}, a compatible $\hat Q_\omega$ can be computed explicitly from $\hat S_\omega$ and any $\hat Q_\omega > 0$ that satisfies \eqref{eq:Q_omega_property}. However, the map $\omega\mapsto\hat Q_\omega$ is generally not guaranteed to be continuous on $\bbR^s$. Instead, we assume continuity and show that, under an additional boundedness condition, $\sys$ is externally equivalent to a well-posed, internally impedance passive and self-dual LTSI system.
	
	\begin{theorem}\label{thm:reciprocity_passivity}
		Suppose that the LTSI system $\sys$ is reciprocal and weakly impedance passive, the LTI system $\hat\sys_\omega$ is minimal for all $\omega \in \bbR^s$, and the maps $\hat S:\bbR^s\mapsto\bbC^{n\times n}$ and $\hat Q:\bbR^s\mapsto\bbC^{n\times n}$ are continuous, pointwise invertible and such that, for all $\omega\in\bbR^s$, $\hat S_\omega = \hat S_\omega\a$ satisfies \eqref{eq:S_omega_property}, $\hat Q_\omega > 0$ satisfies \eqref{eq:Q_omega_property}, and $\hat Q_\omega = \hat S_\omega\hat Q_\omega\inv \hat S_\omega$. If, in addition,
		\begin{equation}\label{eq:CSinvQinvB_bounded}
			\sup_{\omega\in\bbR^s}\norm{\hat C_\omega \hat S_\omega\inv \hat Q_\omega \hat B_\omega} < \infty,
		\end{equation}
		then there exists a well-posed LTSI system
		\begin{equation}\label{eq:ltsi_ss_bar}
			\bar\sys:\left\lbrace\begin{aligned}
				\frac{\d}{\d t}\bar z(t) &= \bar \calA \bar z(t) + \bar \calB u(t),\\
				y(t) &= \bar\calC \bar z(t),
			\end{aligned}\right.
		\end{equation}
		with state space $\bar\calZ = \calZ$ and impulse response $\calG$, that is both internally impedance passive and self-dual. Furthermore, we can partition $\bar\calZ = \bar\calZ_1\times\bar\calZ_2$ and, accordingly,
		\begin{equation}\label{eq:ltsi_ss_bar_partition}
			\bar \calA = \bbm \bar\calA_{11} & \bar\calA_{12} \\ -\bar\calA_{12}\a & \bar\calA_{22} \ebm,\quad \bar\calC\a = \bar\calB = \bbm \bar\calB_1 \\ 0\ebm,
		\end{equation}
		where the operators $\bar\calA_{11}$, $\bar\calA_{12}$ and $\bar\calA_{22}$ are densely defined, $\bar \calA_{11} \leq 0$ and $\bar \calA_{22} \leq 0$.
	\end{theorem}
	\proof
		Since $\hat Q$ is continuous and $\hat Q_\omega>0$ for all $\omega\in\bbR^s$, there exists a pointwise invertible continuous map $\hat L:\bbR^s\to\bbC^{n\times n}$ such that $\hat Q_\omega = \hat L_\omega \hat L_\omega\a$ for all $\omega\in\bbR^s$. (e.g., the pointwise Cholesky factorization). Note that the map $\omega\mapsto \hat L^{-1}_\omega \hat S_\omega \hat L_\omega^{-\ast}$ is continuous on $\bbR^s$, pointwise self-adjoint and invertible. Therefore, there exists a pointwise unitary map $\hat U: \bbR^s\to\bbC^{n\times n}$ and pointwise diagonal and invertible  continuous map $\hat D:\bbR^s\to\bbR^{n\times n}$ such that 
		\begin{equation}
			\hat L^{-1}_\omega \hat S_\omega \hat L_\omega^{-\ast} = \hat U_\omega \hat D_\omega \hat U_\omega\a
		\end{equation}
		for all $\omega\in\bbR^s$. Now, let $\hat T:\bbR^s\to\bbC^{n\times n}$ be the pointwise invertible map given by $\hat T_\omega = \hat U_\omega\a \hat L_\omega\a$ for all $\omega\in\bbR^s$, and note that $\hat Q_\omega = \hat T_\omega\a \hat T_\omega$ and $\hat S_\omega = \hat T_\omega\a\hat D_\omega \hat T_\omega$ for all $\omega\in\bbR^s$. The assumption that $\hat Q_\omega = \hat S_\omega\hat Q_\omega\inv \hat S_\omega$ implies that $I = \hat D_\omega^2$, that is, $\hat D_\omega$ is a signature matrix for all $\omega\in\bbR^s$. Therefore, since $\hat D$ is continuous, it must be constant and, without loss of generality, we can assume that
		\begin{equation}\label{eq:D_omega_signature}
			\hat D_\omega = \bbm I & 0 \\ 0 & -I \ebm
		\end{equation}
		for all $\omega\in\bbR^s$. 
		
		Next, let the maps $\hatbar A:\bbR^s\to\bbC^{n\times n}$, $\hatbar B:\bbR^s\to\bbC^{m\times n}$, and $\hatbar C:\bbR^s\to\bbC^{m\times n}$ be given by
		\begin{equation}\label{eq:hatbar_omega}
			{\hatbar A}_\omega = \hat T_\omega \hat A_\omega \hat T_\omega\inv,\quad \hatbar B_\omega = \hat T_\omega \hat B_\omega,\quad \hatbar C_\omega = \hat C_\omega \hat T_\omega\inv.
		\end{equation}
		for all $\omega\in\bbR^s$. This corresponds to a state-space transformation of $\hat \sys_\omega$ given by $\hatbar z_\omega = \hat T_\omega \hat z_\omega$. It is straightforward to verify that \eqref{eq:S_omega_property} yields
		\begin{equation}
			\hatbar A_\omega\a \hat D_\omega = \hat D_\omega \hatbar A_\omega,\quad \hatbar C_\omega\a = \hat D_\omega\hatbar B_\omega
		\end{equation}
		for all $\omega\in\bbR^s$, while \eqref{eq:Q_omega_property} yields
		\begin{equation}\label{eq:hatbar_internally_passive}
			\hatbar A_\omega\a + \hatbar A_\omega \leq 0,\quad \hatbar C_\omega\a = \hatbar B_\omega
		\end{equation}
		for all $\omega\in\bbR^s$.
		Consequently, due to \eqref{eq:D_omega_signature}, we can write
		\begin{equation}\label{eq:lti_ss_bar_partition}
			\hatbar A_\omega\a = \bbm \hatbar A_{\omega,11} & \hatbar A_{\omega,12} \\ -\hatbar A_{\omega, 12}\a & \hatbar A_{\omega, 22} \ebm,\quad \hatbar C_\omega\a = \hatbar B_\omega = \bbm \hatbar B_{\omega,1} \\ 0 \ebm
		\end{equation}
		where $\hatbar A_{\omega, 11} \leq 0$, $\hatbar A_{\omega,22} \leq 0$.
		
		With this in mind, let $\bar \calA$, $\bar \calB$ and $\bar\calC$ be the translation invariant operators with Fourier symbols $\hatbar A$, $\hatbar C$ and $\hatbar B$. Consider the corresponding LTSI system $\bar\sys$ given by \eqref{eq:ltsi_ss_bar}. To show that $\bar\sys$ is well-posed, we need to show that $\bar \calA:\bar\calZ\supset\D(\bar\calA)\to\bar\calZ$, is the infinitesimal generator of a $C_0$-semigroup, $\bar \calB\in\B(\calU,\bar\calZ)$ and $\bar\calC\in\B(\bar \calZ,\calY)$. From \eqref{eq:hatbar_internally_passive} it follows that $\bar \calA$ is the infinitesimal generator of a contraction $C_0$-semigroup, see \cite[Corollary~2.3.3]{curtain2020}, and $\bar \calB$ is bounded if and only if $\bar \calC$ is bounded. Note that $\hatbar B_\omega = \hat T_\omega \hat B_\omega = \hat T_\omega \hat S_\omega\inv \hat C_\omega\a$, hence 
		\begin{align*}
			\norm{\hatbar B_\omega v}^2 &= \inner{v}{\hat C_\omega \hat S_\omega\inv T_\omega\a  \hat T_\omega \hat B_\omega v} \leq \norm{\hat C_\omega \hat S_\omega\inv \hat Q_\omega \hat B_\omega\a} \norm{v}^2
		\end{align*}
		for all $\omega\in\bbR^s$ and $v\in\bbC^n$. This implies that $\bar\calB$ and, thus, $\bar\calC$, are bounded due to \eqref{eq:CSinvQinvB_bounded}. Now, the well-posed system $\bar\sys$ is internally impedance passive due to \eqref{eq:hatbar_internally_passive}. It is also self-dual with $\bar z\dual = \calD \bar z$ being the unique state-space isomorphism, where $\calD$ is the translation invariant operator with Fourier symbol $\hat D$, which is bounded with bounded inverse due to \eqref{eq:D_omega_signature}. Finally, the partition \eqref{eq:ltsi_ss_bar_partition} follows from \eqref{eq:lti_ss_bar_partition}.
		\endproof
	
	Theorem~\ref{thm:reciprocity_passivity} can be useful even when $\hat \sys_\omega$ is minimal only for $\omega$ in a dense subset of $\bbR^s$. We demonstrate this with the example of the Timoshenko beam.
	\begin{example}
		As shown in Example~\ref{ex:timoshenko_reciprocal}, the (ill-posed) LTSI system $\sys$ that describes the Timoshenko beam is reciprocal. Recall that the total energy is defined via the differential operator $\calQ$ in \eqref{eq:calQ}. The Fourier symbol of $\calQ$ is given by
		\begin{equation}
			\hat Q_\omega = \bbm \omega^2 & -j\omega & 0 & 0 \\  j\omega &  1 +\omega^2 & 0 & 0 \\ 0 & 0 & 1 & 0 \\ 0 & 0 & 0 & 1 \ebm
		\end{equation}
		which is continuous in $\omega$ and satisfies \eqref{eq:Q_omega_property} for all $\omega\in\bbR^s$. Furthermore, we have that $\hat Q_\omega = \hat S_\omega \hat Q_\omega\inv \hat S_\omega$, where $\hat S_\omega$ is given in \eqref{eq:timoshenko_S_omega}. We can verify that $\hat C_\omega\hat S_\omega\inv \hat Q_\omega\hat B_\omega = 1$, hence \eqref{eq:CSinvQinvB_bounded} holds. Note that $\hat Q_\omega = \hat L_\omega \hat L_\omega\a$ with
		\begin{equation}
			\hat L_\omega = \bbm j\omega & 0 & 0 & 0\\ -1 & j\omega & 0 & 0 \\ 0 & 0 & 1 & 0 \\ 0 & 0 & 0 & 1\ebm
		\end{equation}
		which implies that
		\begin{equation}
			\hat L_\omega\inv \hat S_\omega L_\omega^{-\ast} = 
			\bbm -1 & 0 & 0 & 0\\
			0 & -1 & 0 & 0\\
			0 & 0 & 1 & 0\\
			0 & 0 & 0 & 1\ebm,
		\end{equation}
		hence the transformation from the proof of Theorem~\ref{thm:reciprocity_passivity} is given by $\hat T_\omega = \hat L_\omega\a$. Consequently, the translation invariant operator $\calT$ with Fourier symbol $\omega\mapsto\hat T_\omega$ is given by
		\begin{equation}
			\calT = \bbm \frac{\partial}{\partial x} & -1 & 0 & 0\\ 0 & \frac{\partial}{\partial x} & 0 & 0 \\ 0 & 0 & 1 & 0 \\ 0 & 0 & 0 & 1\ebm
		\end{equation}
		hence the state $\bar z = \calT z$ from Theorem~\ref{thm:reciprocity_passivity} is given by \eqref{eq:timoshenko_barz}, i.e., the known physical state of the Timoshenko beam.
	\end{example}
	
	We conclude this section with a few remarks.
	\begin{remark}
		Theorem~\ref{thm:reciprocity_passivity} is simplified in the special cases where $\sys$ is (weakly) lossless or of relaxation type. If $\hat \sys_\omega$ is lossless for all $\omega\in\bbR^s$, then it has a unique storage function, which must be compatible. In fact, this is the case for the (ill-posed) LTSI system describing the Timoshenko beam. On the other hand, if $\hat \sys_\omega$ is of relaxation type for all $\omega\in\bbR^s$, then there is a unique compatible storage and $\hat Q_\omega = \hat  S_\omega > 0$ for all $\omega\in\bbR^s$. This implies that both $\hat S$ and $\hat Q$ in Theorem~\ref{thm:reciprocity_passivity} are continuous because the symbols $\hat A$, $\hat B$ and $\hat C$ are continuous. Moreover, \eqref{eq:CSinvQinvB_bounded} follows from the assumption that $\calB$ and $\calC$ are bounded, and the resulting LTSI system $\bar\sys$ is internally of relaxation type \cite{donchev2025arxiv}, i.e., $\bar\calA \leq 0$ and $\calC\a = \calB$.
	\end{remark}
	
	\begin{remark}
		In view of the partitioning \eqref{eq:ltsi_ss_bar_partition}, the LTSI system $\bar\sys$ in Theorem~\ref{thm:reciprocity_passivity} can be interpreted as an infinite-dimensional port-Hamiltonian system
		\begin{equation}\label{eq:ltsi_ph}
			\bar\sys:\left\lbrace
			\begin{aligned}
				\frac{\d}{\d t}\bar z(t) &= (\calJ - \calR) \bar z(t) + \calG u(t),\\
				y(t) &= \calG\a \bar z(t),
			\end{aligned}
			\right.
		\end{equation}
		with Hamiltonian $\calH(\bar z) = \half\norm{\bar z}^2$, where
		\begin{equation*}
			\calJ = \bbm 0 & \bar\calA_{12} \\ -\bar\calA_{12}\a & 0 \ebm,\quad \calR = -\bbm \bar\calA_{11} & 0 \\ 0 & \bar\calA_{22} \ebm,\quad \calG = \bbm \bar\calB_1 \\ 0 \ebm
		\end{equation*}
		represent the lossless interconnection structure, the dissipative resistive structure, and the port structure, respectively. Note that $\calJ$ and $\calR$ are densely-defined, $\calJ$ is skew-adjoint, and $\calR$ is self-adjoint and nonnegative. This simple extension of (linear) finite-dimensional port-Hamiltonian systems \cite{vanderchaft2014} to the infinite-dimensional setting is not prominent in the literature, which mostly focuses on boundary control systems and specific classes of partial differential equations, see, e.g., \cite{jacob2012, skrepek2021, maschke2023}. Recently, a very general treatment of infinite-dimensional port-Hamiltonian systems that includes \eqref{eq:ltsi_ph} as a special case is presented in \cite{philipp2025}.
	\end{remark}
	\begin{remark}\label{rem:reciprocity_stability}
		The requirement on weak impedance passivity in Theorem~\ref{thm:reciprocity_passivity} can be relaxed to a form of weak Lyapunov stability. Indeed, the proof of Theorem~\ref{thm:reciprocity_passivity} remains essentially the same if we assume that $\hat Q_\omega>0$ satisfies only the left-hand side of \eqref{eq:Q_omega_property}. The resulting LTSI system $\bar\sys$ is still well-posed and self-dual, but not necessarily internally impedance passive. Instead, it is only stable, in the sense that the state operator $\bar\calA$ is the generator of a contraction $C_0$-semigroup.
	\end{remark}
	
	\section{Conclusion}\label{sec:conclusion}
	
	We extended the concept of reciprocity in LTI systems to a class of LTSI systems. Although reciprocity is an external property, we showed that it also constrains the internal structure of the system. Most importantly, we demonstrated that the combination of reciprocity and (weak) impedance passivity enables the construction of a well-posed, self-dual, and internally impedance-passive state-space realization of an otherwise ill-posed LTSI system. This highlights the fact that incorporating external physical properties into state-space modelling is essential, especially in the context of infinite-dimensional systems. 
	
	There are several directions for future research. One question is whether reciprocity, together with other forms of (weak) stability, is enough to guarantee the existence of a well-posed state-space realization of a given LTSI system, see Remark~\ref{rem:reciprocity_stability}. Another is to investigate the relationship between reciprocal and pseudo-gradient LTSI systems. Finally, it would be worthwhile to extend these results to other spatial domains, such as intervals with periodic boundary conditions.
	
	\bibliographystyle{ieeetr}
	\bibliography{../../references/kernel_refs}
\end{document}